  \documentclass{ajour}
 \usepackage{ajourps}

\usepackage{amsmath}

\input{boxedeps.tex}
\SetEPSFDirectory{./}
\HideDisplacementBoxes

\SetRokickiEPSFSpecial  

\DeclareMathAlphabet{\mathgothic}{U}{euf}{m}{n}
\def\gB{\mathgothic{B}}
\def\lcomp{[\kern -1.5pt [}
\def\rcomp{]\kern -1.4pt ]}
\def\s{{\bf s}}

\begin{document}

%



\authorrunninghead{Brent Everitt}
\titlerunninghead{Alternating Quotients of Fuchsian Groups}





\title{Alternating Quotients of Fuchsian Groups}


\author{Brent Everitt
\thanks{Part of this work was done while 
the author was aguest of Sonderforschungsbereich 343,
Unversit\"{a}t Bielefeld. He is grateful for their financial support and hospitality.}}

\affil{Department of Mathematical Sciences, University of Durham,\\ Durham DH1 3LE, 
England}

\email{brent.everitt@durham.ac.uk}

\abstract{It is shown that any finitely generated, non-elementary Fuchsian group has 
among its homomorphic images all but finitely 
many of the alternating groups $A_n$.
This settles in the affirmative a long-standing conjecture of  Graham Higman. }

\begin{article}
\section{Introduction}\label{intro}

It all started with a theorem of G. A. Miller \cite{Miller01}:
the classical modular group $\text{PSL}_2(\mathbf{Z})$ has among its
homomorphic images every alternating group, except $A_6,A_7$ and $A_8$.
In the late 1960's Graham Higman conjectured that any (finitely generated non-elementary) 
Fuchsian group has among its
homomorphic images all but finitely many of the alternating groups.  
This reduces to an investigation
of the cocompact $(p,q,r)$-triangle groups, and in the series of papers
\cite{Conder81,Conder80,Everitt94,Mushtaq93,Rota92} the conjecture was 
verified in the affirmative when $p=2$.
Assuming the Fuchsian group is finitely generated and non-elementary, and
taking the phrases ``almost all" to 
be synonymous with ``all but finitely many", and ``surjects'' with ``has among its homomorphic 
images'', we build on this earlier work to prove 

\begin{proclaim}{Theorem}\label{th1}
Any Fuchsian group surjects almost all  of the alternating groups.
\end{proclaim}

There are several motivations behind the conjecture: Fuchsian groups have an
algebraic structure that is somewhat complicated, and 
to get a firmer grip on this situation, one may be tempted to 
consider their finite, or even simple, homomorphic images. 
There is also a geometric incentive, namely, any compact Riemann
surface (or complex algebraic curve) of genus $>1$ has conformal automorphism group a finite
homomorphic image of some Fuchsian group. 

Schreier coset diagrams supply the technology used to prove the theorem, and
they appear in the literature in various guises (see \cite{Cohen94,Jones96} for alternative
formulations as hypermaps or {\em dessin d'enfants\/}). Section 
\ref{coset.diagrams} has the definition and the basic properties. 
Section \ref{triangle.groups} contains the proof of the theorem.


\section{The plan}\label{genus_nonzero}

Suppose $X$ is the 2-sphere $S^2$, the Euclidean plane $\mathbf{E}^2$ or the hyperbolic plane $\mathbf{H}^2$.
Let $G$
be a finitely generated non-elementary discrete group of orientation preserving isometries of 
$X$. By classical work of Fricke and
Klein (see for instance \cite{Zieschang80}), $G$ has a presentation of the form,
\begin{align*}\label{presentation}
\text{generators:\hspace{1em}}&a_{1},b_{1},\ldots,a_{g},b_{g},&&\text{(hyperbolic)}\\
&x_{1},\ldots,x_{e},&&\text{(elliptic)}\\
&y_{1},\ldots,y_{s},&&\text{(parabolic)}\\
&z_{1},\ldots,z_{t}.&&\text{(hyperbolic boundary elements)}\\
\text{relations:\hspace{1em}}&x_{1}^{m_{1}}=\cdots=x_{e}^{m_{e}}=1,&\\
&\prod_{i=1}^{e}x_{i}\prod_{j=1}^{s}y_j
\prod_{k=1}^{t}z_k
\prod_{l=1}^{g}[a_{l},b_{l}]=1.&
\end{align*}
When $X=\mathbf{H}^2$, $G$ is called a Fuchsian group.
The
division into spherical, Euclidean and Fuchsian is governed by the quantity,
\begin{equation}\label{riemann.hurwitz}
\mu(G)=2g-2+\sum_{i=1}^{e}\biggl(1-\frac{1}{m_i}\biggr)+s+t,
\end{equation}
with $\mu(G)<0, =0$ or $>0$ as $X=S^2,\mathbf{E}^2$ or $\mathbf{H}^2$.
The quotient $X/G$ is an orientable 2-orbifold of genus $g$ with $e$ cone points, $s$ punctures
and $t$ boundary components.  Its geometry and the algebraic structure of $G$ are intimately
connected, so that $G$ is determined upto isomorphism by its signature
$(g;m_1,\ldots,m_e;s;t)$, $2\leq m_1\leq\cdots \leq m_e$.

To prove the theorem, it suffices to just consider the cocompact Dyck groups--the
cases where in the signature we have $g=s=t=0$. To see why we make a few elementary 
observations.
\begin{enumerate}
\item A group of signature $(g;m_1,\ldots,m_e;s;t)$ is isomorphic to one of
$(g;m_1,\ldots,$ $m_e;s+t;0)$, and by (\ref{riemann.hurwitz}), the former is Fuchsian if
and only if the latter is. We may assume then that $t=0$. Write
$(g;m_1,\ldots,m_e;s)$ instead of 
$(g;m_1,\ldots,m_e;s;0)$ from now on.
\item We can surject $G=(g;m_1,\ldots,m_e;s)$ onto 
$G'=(g';m_1,\ldots,m_i',$ $\ldots,\hat{m_j},$ $\ldots,m_e;s')$, for any $g'\leq g$, 
$s'\leq s$, and $m_i'$ a
divisor of $m_i$. The hat denotes ommission. Here's how: 
map the $j$-th elliptic, $s-s'$ of the parabolic and $g-g'$ hyperbolic pairs
of generators of $G$ to the identity of $G'$; map the $i$-th elliptic generator of $G$
to the corresponding elliptic generator of $G'$ raised to the power $m_i/m_i'$.
All other generators of $G$ map to the corresponding ones in $G'$.
The map then extends to the desired homomorphism.
\item Writing $(m_1,\ldots,m_e)$ when $g=s=0$, suppose $\psi:G=(m_1,\ldots,m_e)
\rightarrow S_n$ is a homomorphism with transitive image and let $G_1$ be the subgroup
of $G$ consisting of those elements stabilising some fixed point of $\{1,2,\ldots,n\}$.
By theorem 1 of \cite{Singerman70}, $G_1$ has signature
$(g';n_{11},n_{12},\ldots,n_{1\rho_1},\ldots,n_{r1},n_{r2},\ldots,n_{r\rho_r})$,
where $\psi(x_i)$ has exactly one cycle each of lengths 
$m_i/n_{i1},\ldots,m_i/n_{i\rho_i}$, with all other cycles of length $m_i$,
and $\mu(G_1)=n\mu(G)$. Moreover, if $G_1$ is normal in $G$, and we have the 
theorem for $G$, the simplicity of $A_n$ for $n\geq 5$ gives the result for
$G_1$ as well.
\item Finally, any $k$-cycle $(a_1,\ldots,a_k)\in A_n$ can be written as a product 
$$(a_2,a_k)(a_3,a_{k-1})\ldots(a_{k/2+1},a_{k/2+2})(a_1,a_2)(a_3,a_k)\ldots(a_{k/2+1},a_{k/2+3}),$$
of two involutions in $A_n$. 
Similarly any cycle of even length in $S_n$ can be written as a product of a involution in $S_n$
and an involution in $A_n$.
Thus, if we have the result for $(m_1,\ldots,k,$ $\ldots,m_e;s)$ we have it for
$(m_1,\ldots,2,2,\hat{k},\ldots,m_e;s)$ too.
\end{enumerate}

\begin{lemma}\label{prop.all}
The theorem is true for every Fuchsian group if it holds for every Dyck group.
\end{lemma}

\begin{proof}
Proceeding according to the genus, suppose $G$ has signature $(g,m_1,\ldots,$ $m_e;s)$ with
$g\geq 2$. Map $G$ onto $\langle x,y\,|\text{---}\rangle$, free of rank two, by sending $a_1\mapsto
x$, $a_2\mapsto y$, and all the other generators to the identity. Since $A_n$ is $2$-generated for
$n\geq 3$ (see \cite{Coxeter65}), we are done.

A group of genus one with $e\geq 1$ can be surjected onto 
$(1;m_1;0)$ for $m_1\geq 2$, by comment 2 above. The map $\theta:(0;2,2,2,2m_1;0)\rightarrow S_2$ sending all generators to
the permutation $(1,2)$ has kernel isomorphic to $(1;m_1;0)$ by
comment 3 above, hence the result holds for groups of genus one with $e\geq 1$.   
For groups of genus one with no periods,
hence signature $(1;\text{---};s)$ for $s\geq 1$, we may surject onto $(1;\text{---};1)$. But this is
easily seen to be free of rank
two, so the result holds here also.

A group of genus zero with no periods must, by (\ref{riemann.hurwitz}),
have at least three parabolic generators, and hence surject $(0;\text{---};3)$.
But this is free of rank two also. With a single period we have $s\geq 2$, and the group surjects
$(0;m_1;2)\cong\mathbf{Z}_{m_1}*\mathbf{Z}$, the free product of $\mathbf{Z}_{m_1}$
and $\mathbf{Z}$. This surjects $\mathbf{Z}_{m_1}*\mathbf{Z}_3$, which in turn surjects any
Fuchsian triangle group of the form $(0;3,m_1,r;0)$.

With two periods and one parabolic, we have $(0;m_1,m_2;1)\cong\mathbf{Z}_{m_1}*\mathbf{Z}_{m_2}$, where $m_2\geq 3$,
so we can surject any Fuchsian triangle group like $(0;m_1,m_2,r;0)$. A group with more parabolics,
$(0;m_1,m_2;s)$ for $s\geq 2$, surjects $(0;m_1;2)$ done above. Finally,
$(0;m_1,\ldots,m_e;s)$, $e\geq 3$, surjects either $\mathbf{Z}_2*\mathbf{Z}_2*\mathbf{Z}_2$ or 
$\mathbf{Z}_{m_1}*\mathbf{Z}_{m_2}$ for $m_2\geq 3$. Surject the former onto a Fuchsian
$(0;2,2,2,p;0)$. The latter has already been handled.
\end{proof}

\begin{lemma}\label{prop_dyck}
The theorem holds for every Dyck group if it holds for the following:
\begin{enumerate}
\item The Fuchsian triangle groups $(p,q,r)$ with $2\leq p<q<r$ distinct primes;
\item the triangle groups $(2,4,r)$ for $r\geq 5$ a prime;
\item the groups $(2,3,8)$, $(2,3,9)$, $(2,3,10)$, $(2,3,12)$, $(2,3,15)$, $(2,3,25)$, $(2,4,6)$, 
$(2,4,8)$, $(2,4,9)$, $(2,5,6)$, $(2,5,9)$ and $(3,4,5)$;
\item the groups $(2,3,3,3)$, and $(3,3,3,3)$.
\end{enumerate}
\end{lemma}

\begin{proof}
The hyperbolic triangle group $(2,m_1,m_2)$ surjects
$(2,q,r)$ for $q$ and $r$ some prime divisors of $m_1$ and $m_2$. 
If $(2,q,r)$ is Fuchsian, we have by (\ref{riemann.hurwitz}) that $1/q+1/r<1/2$.
If $q$ and $r$ are distinct, we have a group listed
in part 1 of the lemma. If $q=r$, the map $\psi:(2,q,4)\rightarrow S_2$ 
that sends the generators
of orders $2$ and $4$ to the permutation $(1,2)$ and the generator of order $q$ to the identity has
kernel $(q,q,2)\cong(2,q,q)$. We have $2/q<1/2$, hence $q\geq 5$, and the 
theorem holds for $(2,q,q)$ if it holds for $(2,4,q)$, a group listed in part 2 of the 
lemma.
 
If $(2,q,r)$ isn't Fuchsian, it must be, after a possible reordering, one of $(2,2,r)$ for
$r\geq 2$, $(2,3,3)$ or $(2,3,5)$.  
The first gives that $(2,m_1,m_2)$ must have the form $(2,m_1,2^{l})$, 
for $m_1\geq 3$ and $l\geq 2$. 
If $m_1=3$ or $4$ then $l\geq 3$, as $(2,3,4)$ is spherical and $(2,4,4)$ Euclidean, so
the group surjects $(2,3,8)$ or $(2,4,8)$, both of which are listed in 
the lemma. For $m_1\geq 5$, $(2,m_1,2^{l})$ surjects
$(2,m_1,4)\cong(2,4,m_1)$.
This in turn surjects $(2,4,r)$, $r$ prime, and we have a 
group listed in part 2 unless
$r=2$ or $3$. In the first case, $m_1=2^{n}\geq 8$, so $(2,4,m_1)$ 
surjects $(2,4,8)$. In the second,
$m_1=2^{l_1}3^{n_1}$, and the group surjects $(2,4,9)$ when 
$l_1=0$, or $(2,4,6)$ otherwise. The cases $(2,q,r)=(2,3,3)$ or $(2,3,5)$ are entirely
similar.

This accounts for the $(2,m_1,m_2)$ Fuchsian groups, and the case of a general triangle
groups is much the same. Similarly for the groups with four or five elliptic 
generators--either they can be surjected directly onto triangle groups or eliminated
from consideration using comment 4 at the beginning of the section. The only exceptions
are those listed in the lemma. Finally, a group with six or more elliptic 
generators can always be surjected directly onto a Fuchsian group with five. No 
doubt the reader can fill in the details.
\end{proof}

In \cite{Conder81,Conder80,Everitt94}, the groups $(2,3,r)$ for all $r\geq 7$
and $(2,4,r)$ for all $r\geq 5$ were dealt with. Theorems 1-3 of \cite{Mushtaq93} take care of 
the $(2,q,r)$,
$5\leq q<r$ prime, with the exception of sixty cases. These sixty, and those
from parts 3 and 4 of Lemma \ref{prop_dyck} can be found in the preprint version
of this paper \cite[$\S 6$]{Everitt98}. This
leaves the triangle groups $(p,q,r)$, $3\leq p<q<r$ to consider, and they can be found in Section 
\ref{triangle.groups}.

Later on we will construct permutation groups as homomorphic images of
Fuchsian groups and will identify the images as 
alternating using,

\begin{theorem}[\cite{Jordan73}, refer to \cite{Wielandt64} Theorem 13.9]
Let $G$ be a primitive permutation group of degree $n$ 
containing a prime cycle for some prime $q\leq n-3$. Then $G$ is
either the alternating group $A_{n}$ or the symmetric group $S_{n}$. 
\label{theorem.jordan}\end{theorem} 

The following lemma, well known to the cognoscenti, allows one to replace primitivity
by more easily verifiable criteria. Recall that the support of a permutation 
$\sigma\in S_n$ consists of those elements of
$\{1,2,\ldots,n\}$ not fixed by $\sigma$.

\begin{lemma}\label{primitive}
Let $G=\langle \sigma_1,\sigma_2,\ldots,\sigma_k\rangle$ be a transitive permutation group
of degree $n$ containing a prime cycle $\mu$. For each $\sigma_i$, suppose there is
a point in the support of $\mu$ whose image under
$\sigma_i$ is also in the support of $\mu$. Then $G$ is primitive. 
\end{lemma}

\begin{proof}
Suppose on the contrary that $G$ is imprimitive with block system $\gB$. For $\sigma\in G$, let
$\bar\sigma$ be the permutation induced by $\sigma$ on $\gB$, and $\overline G$ the group generated by
the $\bar\sigma_i$. The map $\sigma\mapsto \bar\sigma$ is an
epimorphism from $G$ onto $\overline G$, and $\overline G$ acts transitively on $\gB$. All blocks
$B\in\gB$ thus have the same size, say $|B|$. If $B\in\gB$ is in $\mbox{supp}(\bar\mu)$, the
support of $\bar\mu$, then $B$ and its image under $\mu$ are distinct blocks, and so $B$ is
contained in $\mbox{supp}(\mu)$.  Taking the union of all the blocks in
$\mbox{supp}(\bar\mu)$ thus gives \begin{equation}\label{supp}
|B||\mbox{supp}(\bar\mu)|\leq |\mbox{supp}(\mu)|.
\end{equation}
Now $\mu$ has order $q$ a prime, and $\bar\mu$ is a homomorphic image of $\mu$. Thus, if
$\bar\mu\not= 1$, then $\bar\mu$ has order $q$, and so $|\mbox{supp}(\bar\mu)|\geq q$. Since
$\gB$ is non-trivial, we have $|B|>1$, and hence by (\ref{supp}),
$|\mbox{supp}(\mu)|>q$. This contradicts the fact that $\mu$ is a $q$-cycle, so we must have
$\bar\mu=1$. This means that $\bar\mu$ fixes every block, or equivalently, any point and its image
under $\mu$ lie in the same block. But $\mu$ is a single cycle, so there
is a block $B^*$ with $\mbox{supp}(\mu)\subseteq B^*$. By the condition stated in the Lemma,
$B^*$ and its image under $\sigma_i$ intersect for all $i$, so are equal. Since the $\sigma_i$
generate $G$, the whole group must fix $B^*$, and by transitivity, $B^*=\{1,2,\ldots,n\}$, so there
is just one block. This is the desired contradiction. \end{proof}


\section{Coset diagrams}\label{coset.diagrams}

Suppose $G$ is a group with a finite presentation $\langle X;R\,\rangle$, and let $K_0=K_0(X;R)$ be the standard
2-complex with $\pi_1(K_0)\cong G$. 
The 1-skeleton of $K_0$ consists of a single vertex incident with oriented loops or 
{\em edges\/}
that are in one to one correspondence with the generators $X$. Each edge $x\in X$ is a pair of
oppositely oriented  arcs, an $x$-arc and an $x^{-1}$-arc. The former coincides with the edge under
its given orientation and the latter to the edge with the reverse orientation. The faces of $K_0$
are in one to one correspondence  with the relators $R$, and are obtained by sewing discs onto the
1-skeleton, each with boundary label a  relator word $r\in R$,
see \cite[\S6.3]{Hilton62}.

A Schreier coset diagram for $G$ is a cellular (that is, $k$-cells lift to $k$-cells)
covering of $K_0$ (see
\cite[\S2.2.1 and \S4.3.2]{Stillwell80} or \cite{Cohen89}). A covering $K$ realises a subgroup
$H\cong \pi_1(K)$ of
$G$, with the vertices of $K$ in one to one correspondence with the cosets of $H$ in $G$. 
Conversely, every subgroup is
realisable in this way from some diagram. 

Their usefulness for our purposes stems from the fact that any coset diagram $K$ yields a homomorphism
$\theta_{K}:G\rightarrow\mbox{Sym}\{\text{vertices of $K$}\}\cong S_n$. Here $n$ is the sheet number of the 
covering, hence the number $|K|$ of vertices in $K$. 
For any $g\in G$ the image of vertex $v$ under the permutation $\theta_{K}(g)$ is the terminal vertex of the 
path starting at $v$ with label $g$. In particular, $\theta_{K}(G)$ is transitive if and only
if $K$ is path-connected.

All of which is, of course, well known. The simplicial complexes that form coset diagrams for $G$
are characterised by two simple properties: 
\begin{enumerate}
\item For each vertex $v$ and generator $x\in X$, there is precisely one $x$-arc and one
$x^{-1}$-arc having initial vertex $v$.
\item The boundaries of the faces are precisely the paths obtained by starting at some vertex
$v$ and traversing a path with label some $r\in R$.
\end{enumerate}

Condition 2 indicates that in their unrefined form, coset diagrams will be a little 
unwieldy--there will be many faces sharing the same set of boundary edges.
To alleviate matters, we use an equivalent construct, suggested by Higman and
used in \cite{Conder81,Conder80,Everitt97,Everitt94,Mushtaq93,Rota92}. It is what
results by identifying such multiple faces.

\begin{figure}
\begin{center}
\begin{picture}(50,80)
\put(23,10){$v$}
\put(77,47){$x_{e-1}$}
\put(77,-7){$x_{e-1}$}
\put(48,74){$x_1$}
\put(-4,74){$x_1$}
\put(-37,47){$x_2$}
\put(-37,-7){$x_2$}
\put(25,20){\circle*{6}}
\qbezier[35](-1.84,6.58)(25,-15)(51.84,6.58)
\put(25,20){\line(2,1){25}}
\put(69.72,42.36){\vector(-2,-1){25}}
\put(25,20){\vector(-2,1){25}}
\put(2.64,31.18){\line(-2,1){25}}
\put(25,20){\vector(2,-1){25}}
\put(47.36,8.82){\line(2,-1){25}}
\put(25,20){\line(-2,-1){25}}
\put(-19.72,-2.36){\vector(2,1){25}}

\put(25,20){\vector(1,2){13.5}}
\put(38,46){\line(1,2){11.5}}
\put(0,70){\vector(1,-2){13}}
\put(25,20){\line(-1,2){12.5}}
\end{picture}
\caption{}\label{rotation}
\end{center}
\end{figure}

Let $G=(m_1,\ldots,m_{e})$ be some fixed but arbitrarily
chosen Dyck group.  A more convenient presentation than given in the introduction is,
$$
\langle x_1,x_2,\ldots,x_{e-1}\,|\,x_1^{m_1}=x_2^{m_2}=\cdots=x_{e-1}^{m_{e-1}}=
(x_1x_2\ldots x_{e-1})^{m_{e}}=1 \rangle.
$$
A $G$-graph is a directed graph with edges labelled $x_1,\ldots,x_{e-1}$ satisfying property (1)
above. Ordering the edges incident with every
vertex as shown in Figure \ref{rotation} yields a 2-cell embedding of a $G$-graph into a closed orientable
surface (see  \cite{White73} for more details on graph embeddings). Each face of this surface
complex $S$ will have boundary label some power of $x_i$ or $x_1x_2\ldots x_{e-1}$. Call $S$ a  
{\em $G$-diagram\/} if for each face, this power divides the order of the appropriate word given
in the presentation. 

In a $G$-diagram, a path starting at $v$ with label $x_i^{m_i}$ or $(x_1\ldots
x_{e-1})^{m_e}$ circumnavigates a face an integral number of times.
Taking the
underlying $G$-graph and sewing in a 2-cell for each such vertex--relator pair yields 
a coset diagram for $G$. Conversely, the  1-skeleton of a coset diagram  is
a $G$-graph in which a path from any vertex with label a relator is closed (as it bounds a face). 
Embedding the graph as above gives a $G$-diagram. We therefore have

\begin{lemma}
A coset diagram for $G$ yields a unique $G$-diagram, and vice-versa. 
\end{lemma}

Consequently, we use the same terminology for $G$-diagrams as for coset diagrams.
In particular, call a face an $x_i$-face or $(x_1\ldots x_{e-1})$-face whenever it has boundary label some power of $x_i$ or $x_1\ldots x_{e-1}$.

The key property of $G$-diagrams, as Higman observed, is that they can sometimes be 
combined to form new ones.
For this we use {\em handles\/}, that is, pairs of vertices
$\alpha$ and $\beta$, each incident with $x_1$-loops, 
so that the path starting at $\alpha$ with label $x_1\ldots x_{e-1}$ terminates at $\beta$. 

Let $K_1,\ldots,K_t$, $t\leq m_1$, be a collection of 
disjoint $G$-diagrams, and the $2m_1$ distinct vertices $\alpha_1,\beta_1,
\ldots,$ $\alpha_{m_1},\beta_{m_1}$ a collection of $m_1$ handles with at least one 
in each diagram.
Take the disjoint
union of all the underlying $G$-graphs, remove the $x_1$-loops at the vertices $\alpha_j$ and $\beta_j$,
and replace them by $x_1$-edges from $\alpha_j$ to $\alpha_{j+1}$ and $\beta_j$ to $\beta_{j-1}$,
subscripts taken modulo $m_1$. Embed the graph in the usual way, and call the resulting surface complex 
$\lcomp K_1,\ldots, K_t\rcomp$ the {\em composition\/}
of $K_1,\ldots,K_t$.

\begin{proposition}\label{comp_prop}
$\lcomp K_1,\ldots, K_t\rcomp$ is also a $G$-diagram with $\sum|K_i|$ vertices.
\end{proposition}

\begin{proof}
The underlying graph of $\lcomp K_1,\ldots, K_t\rcomp$ is clearly a $G$-graph, so it remains to
show that all faces have boundary labels of the required form. If the boundary of a face does not
contain an
$x_1$-edge with initial vertex one of the $\alpha_j$ or $\beta_j$, then all edges are contained in a
single $G$-diagram $K_i$, and we are done.

Otherwise, we obtain the boundary label for the face by starting at
an $\alpha_j$ or $\beta_j$ and traversing a path with label some power of $x_1$ or some power of 
$x_1\ldots x_e$, until it closes (which it does by repeating an arc). The path obtained by
traversing just $x_1$-edges passes through the vertices
$\alpha_{j+1},\ldots,\alpha_{m_1},\alpha_1,\ldots,\alpha_j$ or 
$\beta_{j-1},\ldots,\beta_1,\beta_{m_1},\ldots,\beta_j$,  
before closing with label $x_1^{m_1}$, so such faces are as they should be.
Observe that {\em before\/} composition, the path starting at $\alpha_j$ 
with label some power
of $x_1\ldots x_e$ arrived at vertex $\beta_j$ after $e$ directed edges, proceeded to traverse
the $x_1$-loop at $\beta_j$ and then an $x_2$-edge. After composition, the path  
from $\alpha_j$ with such a label
arrives instead at $\beta_{j+1}$ after $e$ directed edges, traverses the new $x_1$-edge
to $\beta_j$, and is then identical with the path before composition.
So the boundary label behaves as if the composition never happened, and is thus of the required form.
The number of vertices is obvious.
\end{proof}

Now suppose $G$ is the triangle group
$$
\langle x,y \,|\, x^p=y^q=(xy)^r=1\rangle,\,\,\, 3\leq p<q<r,
$$
with $p,q$ and $r$ prime.
In practice, we 
simplify $(p,q,r)$-diagrams when drawing them: a shaded $q$-gon indicates a $y$-face with
boundary label $y^q$, and a shaded wedge \BoxedEPSF{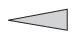 scaled 650} a $y$-face with label
$y$; the orientation on arcs runs anticlockwise around any face they bound unless indicated otherwise; 
$x$-faces with boundary $x$ are removed completely, leaving only the incident vertex which will be called {\em free\/}. On occasion, we will
talk of attaching $x$-arcs to free vertices, by which we mean attach the arcs to the
underlying $G$-graph and re-embed.  

As a consequence, the unshaded faces are precisely the $x$ and $xy$-faces,
and for an embedded $G$-graph to 
be a $G$-diagram, it is sufficient that the $xy$-faces have a number of $y$-arcs dividing 
the appropriate order in their
boundaries, and the $x$-faces a number of $x$-arcs similarly. 
These criteria can usually 
be verified at a glance. 

\begin{figure}
\begin{center}
\BoxedEPSF{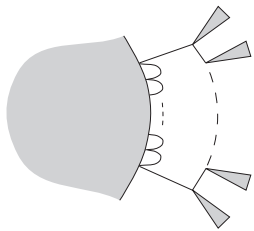 scaled 1000}
\caption{\hspace{0.5em}{\rm Type $k$ pendant.}}
\label{pendant1}
\end{center}
\end{figure}

We devote the remainder of this section to diagrams for
triangle groups. An $x$-face is of type
$[l_1,\ldots,l_{\lambda},\ldots,\bar l_{\mu},\ldots,l_t]$, $\sum_{i=1}^{t} l_i=p$, if it has
boundary label $x^p$, and in traversing the boundary with the orientation,  
\begin{itemize}
\item vertices $(\sum_{i<\lambda} l_i +1)$ through $(\sum_{i\leq \lambda} l_i)$ are
consecutive on some $q$-gon, 
\item vertices $(\sum_{i<\mu} l_i +1)$ through $(\sum_{i\leq \mu} l_i)$ are incident with 
\BoxedEPSF{222p.fig7.eps scaled 650}'s.
\end{itemize}
Of course the face also has type $X$ for $X$ any cyclic permutation of the $l_i$, but in 
practice this
ambiguity causes no confusion. We tend to say {\em type $X$ $x$-cycle\/} rather than $x$-face of type
$X$. Figure
\ref{pendant1} shows a type
$[k,\overline{p-k}]$ $x$-cycle,
$1\leq k\leq p$, or {\em type $k$ pendant\/}.
 
Suppose we have $k$ consecutive free vertices on a $q$-gon, all in the
boundary of the same $xy$-face $F$. Attaching a type $k$ pendant to these vertices
increases the number of $y$-arcs in the boundary of $F$ by $p-2k+1$. The modification also
produces a new $x$-face with boundary $x^p$ and some $y$ and $xy$-faces with
label $y$ and $xy$.

Suppose $q=lp+s$ for $l\geq 1$ and $1\leq s\leq p-1$. Take a shaded $q$-gon, and attach $l-1$ type $p$ pendants
to $p(l-1)$ consecutive vertices. Attach a single type $s$ pendant so that $p$
consecutive vertices are left free. The resulting $q$-gon together with the attachments will be called a 
{\em booster\/}.

\begin{figure}
\begin{center}
\BoxedEPSF{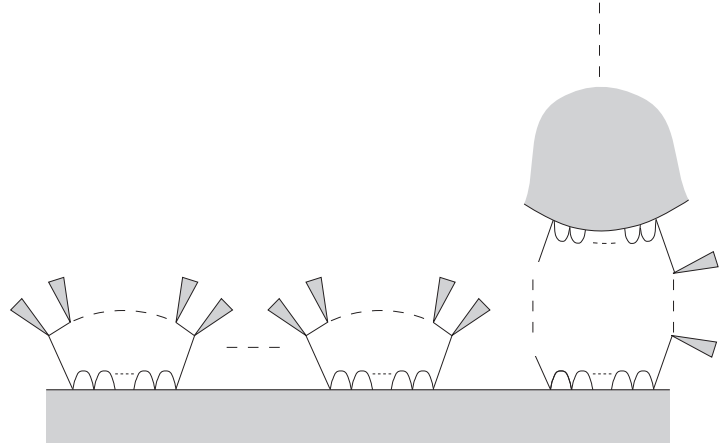 scaled 1000}
\caption{\hspace{0.5em}{\rm Type $\{k_1,\ldots,k_t;X_1,\ldots,X_m\}$ array.}}
\label{array}
\end{center}
\end{figure}

Let
$X_i=[l_{i1},\bar{l}_{i2},{l}_{i3},\bar{l}_{i4}]$, $i=1,\ldots,m$. Suppose that for integers 
$1\leq k_1,\ldots,k_t\leq p$, we have $l_{11}+\sum k_i$ consecutive
free vertices on a $q$-gon bounding an $xy$-face $F$. By attaching a
type $\{k_1,\ldots,k_t;X_1,\ldots,X_m\}$ array to these free
vertices we mean,
\begin{itemize}
\item attach $t$ pendants of types $k_1,\ldots,k_t$, and
\item a collection of $m$ boosters, joined into a chain, with $l_{i3}$ vertices of the $i$-th 
booster connected to
$l_{i1}$ vertices of the $(i-1)$-st by an $x$-cycle of type $X_i$ (taking the $0$-th booster
to be the original $q$-gon)-see Figure \ref{array}.
\end{itemize}
Write $\{k_1,\ldots,k_i^{\delta_i},\ldots,k_t;X_1,\ldots,X_j^{\delta_j},\ldots,X_m\}$ when
the array includes $\delta_i$ type $k_i$ pendants and $\delta_j$ $x$-cycles of type $X_j$. 
Notice that a type
$\{k;\text{---}\}$ array is merely a type $k$ pendant. In attaching an array, the number of
$y$-arcs in the boundary of $xy$-face $F$ increases by 
\begin{equation}\label{array.number}
m(p+l+2-s)+\sum_{i=1}^t (p-2k_i+1)+2\sum_{\bar{l}_{ij}\in X_i} l_{ij},
\end{equation}
together with the creation of the usual complement of $x,y$ and $xy$-faces having boundary $x,x^p,y$
and $xy$. All other faces are unaffected. To see (\ref{array.number}), start with each
$X_i=[1,p-1]$, and observe that replacing it by $[1,\bar{1},p-2]$ increases the $y$-arc count by
two, while a change to
$[2,p-2]$ has no effect.

If $K$ is a $(p,q,r)$-diagram with $g\in (p,q,r)$, the cycle structure of
$\theta_{K}(g)$ is a function $\s:\mathbf{Z}^+\rightarrow\mathbf{Z}^+\cup\{0\}$, such that
$\s(i)$ is the number of cycles of length $i$  when $\theta_{K}(g)$ is written as a
product of disjoint cycles. 
Given two structures $\s_1$ and $\s_2$,
let $\s_1+\s_2$ be their pointwise sum as functions. In Section \ref{triangle.groups}
we will be interested in the structure of the element $x^{-1}y$.

\begin{lemma}\label{cycle.structure}
Suppose $K_1,\ldots,K_t$ are $(p,q,r)$-diagrams with $\s_i$ the cycle structure of
$\theta_{K_i}(x^{-1}y)$. If $K=\lcomp K_1,\ldots,K_t\rcomp$, then 
$\theta_{K}(x^{-1}y)$ has cycle structure $\sum \s_i$.  \end{lemma}

\begin{proof}
Only cycles in $\theta_{K_i}(x^{-1}y)$ that pass through handle points are affected by
the composition. If $\alpha_j$ and $\beta_j$ lie in such a cycle, then in
$\theta_{K}(x^{-1}y)$ the cycle is identical, except that $\beta_j$ is replaced by $\beta_{j-1}$.\end{proof}

Consequently, consideration of the cycle structure of $\theta_{\lcomp K_1,\ldots,K_t\rcomp}(x^{-1}y)$ 
reduces to an investigation of the $\theta_{K_i}(x^{-1}y)$.

\begin{figure}
\begin{center}
\begin{tabular}{ccc}
\BoxedEPSF{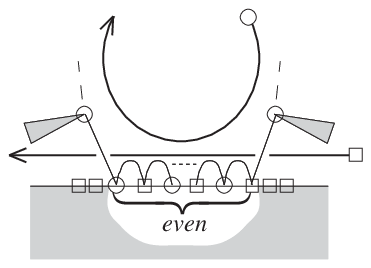 scaled 1000}&
\vrule width 20 mm height 0pt depth 0 pt&\BoxedEPSF{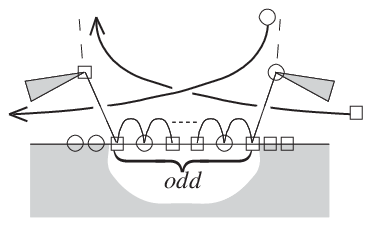 scaled 1000}\\
\end{tabular}
\caption{} \label{cycle.structures}
\end{center}
\end{figure}

We determine the effect on $\theta_{K}(x^{-1}y)$ of attaching an array
by considering the various ingredients. From now on, when we talk of a {\em cycle\/} in
$K$, we will mean a cycle of $\theta_{K}(x^{-1}y)$, and the context should make clear
which cycle we mean. Notice first that consecutive free vertices on a 
$q$-gon are contained in the same cycle. Attaching a type $k$ pendant 
to these vertices increases the length of this cycle by $p-k$ when $k$ is odd. When $k$ is even, it
decreases by $k/2$, and a new cycle of length $p-k/2$ is created. Next, the vertices 
of an isolated
booster are organised into a single cycle of length
$$
q+
\left\{\begin{array}{l}
p-s, s\text{ odd}\\
-{\displaystyle \frac{s}{2}}, s\text{ even.}
\end{array}\right. 
$$
When contained as the $i$-th booster of an array, vertices may be gained or lost from this cycle (it
may even be fused with cycles from  neighbouring boosters) depending
on whether $l_{i1}$ and $l_{i3}$ are even or odd. Figure \ref{cycle.structures} shows the
possible orbits on the vertices, illustrated by small circles and squares.

It will be useful to have at our disposal various maneuvers in which an array is replaced by
another. Replacing an array of type $\{k_1,\ldots,k_t;X_1,\ldots,X_m\}$ by one of type 
$\{k_1,\ldots,k_t,\frac{p+1}{2};$ $X_1,\ldots,X_m\}$ is called spoiling. A push-pull
substitutes $\{k_1,\ldots,$ $k_i-1,\ldots,k_j+1,\ldots,k_t;
X_1,\ldots,$ $X_m\}$, while replacing by 
$\{k_1,\ldots,k_t;X_1,\ldots,X'_i,\ldots,X_m\}$, where $X'_i=
[l_{i1}\pm 1,\bar{l}_{i2},l_{i3}\mp 1,
\bar{l}_{i4}]$, will be known as modifying a chain. 

A few brief remarks on each then. Suppose $K'$ is the result of performing such a
maneuver on some array in the $(p,q,r)$-diagram $K$:
\begin{itemize}
\item $K\xrightarrow{\text{spoiling}}K'$: since (\ref{array.number}) is unchanged, $K'$ is
also a $(p,q,r)$-diagram. The modification requires $\frac{p+1}{2}$ free vertices and
$|K'|= |K|+\frac{p-1}{2}$. The length of the cycle containing these free vertices changes
by a non-trivial amount $<q$.
\item $K\xrightarrow{\text{push-pull}}K'$: again (\ref{array.number}) is invariant so $K'$
is a $(p,q,r)$-diagram. No free vertices are required and $|K'|=|K|$. The length of 
the cycle on the $q$-gon 
to which the array is attached changes by 
\begin{equation}\label{push.pull}
\sum_{\substack{k\in\{k_i,k_j\}\\\text{$k$ even}}} \biggl(p-\frac{k}{2}\biggr)
-\sum_{\substack{k\in\{k_i,k_j\}\\\text{$k$ odd}}} \frac{k}{2}.
\end{equation}
\item $K\xrightarrow{\text{modifying chain}}K'$: again $K'$ is a $(p,q,r)$-diagram, 
with $|K'|=|K|$.
The operation requires a free vertex on the $(i\mp 1)$-st booster, creating one 
on the $(i\pm 1)$-st. Use Figure \ref{cycle.structures} to monitor the effect on cycles in
$\theta_{K}(x^{-1}y)$.
\end{itemize} 


\section{The proof of the theorem}\label{triangle.groups}

Higman's construction, forming the basis of \cite{Conder81,Conder80,Everitt97,Everitt94,Mushtaq93,Rota92}, 
is essentially,

\begin{proposition}
Let $K_1,K_2$ and $K_3$ be path-connected diagrams for the triangle group $(p,q,r)$
such that,
\begin{enumerate}
\item $|K_1|, |K_2|$ are relatively prime, and $|K_3|\geq q+3$;  
\item $K_1$ and $K_2$ each contain at least two handles and $K_3$ one;
\item if $\s_i$ is the cycle structure of $\theta_{K_i}(x^{-1}y)$, then 
$\s_1(kq)=\s_2(kq)=0$, $k\geq 1$, and 
$$\s_3(kq)=
\left\{\begin{array}{ll}
1,&k=1,\\
0,&k>1;
\end{array}\right.
$$
\item if $\mu$ is the $q$-cycle in $\theta_{K_3}(x^{-1}y)$ there are $i,j\in\mu$, not
contained in the handle, with $i^x,j^y\in\mu$.
\end{enumerate}
Then $G=(p,q,r)$ surjects almost all of the alternating groups. 
\end{proposition}

\begin{proof}
Let $p_1,p_2>p$ be distinct primes 
not dividing $|K_1|$ and $|K_2|$. For $k_1$ and $k_2$ arbitrary non-negative integers we construct 
a sequence of diagrams $C_0,C_1,\ldots,C_{k_1},\ldots,C_{k_1+k_2}:=K$ as follows:
for the $0$-th step, if either $k_1\text{ or }k_2=0$, take $C_0=K_3$, otherwise, $C_0=K_2$. 
At step $i$, $1\leq i\leq k_1$, take $p_1$ identical copies of $K_1$ and let $C_i$ be the 
composition,
\begin{equation}\label{composition}
\lcomp\lcomp\ldots\lcomp\lcomp C_{i-1},\underbrace{K_1,\ldots,K_1}_{p-1}\rcomp,
\ldots\rcomp,
\underbrace{K_1,\ldots,K_1}_{p-1}\rcomp,
\underbrace{K_1,\ldots,K_1}_{\leq p-2}\rcomp.
\end{equation}
In particular, the two handles on each $K_1$
allow us to perform the composition, which is a $(p,q,r)$-diagram by Proposition \ref{comp_prop}.
Observe that $C_i$ has at least two handles.
At step $i$,
$k_1+1\leq i\leq k_1+k_2-1$,
take $p_2$ identical copies of $K_2$ and let $C_i$ be a composite diagram 
of the form (\ref{composition}) but with $p_2$ 
copies of $K_2$
instead of $p_1$ copies of $K_1$.
Finally, at step $k_1+k_2$, if $k_1\text{ or }k_2=0$, let $C_{k_1+k_2}$ be as in the previous step. Otherwise, take
a diagram of the form (\ref{composition}) but replace one of the $K_2$'s by a $K_3$ (using its
sole handle).

A quick sketch may help the reader to see what is going on.
Now $|K|=k_1p_1|K_1|+k_2p_2|K_2|+|K_3|$, and
since $|K_1|$ and $|K_2|$ are relatively prime, so too are 
$p_1|K_1|$ and $p_2|K_2|$. By choosing $k_1$ and $k_2$ suitably,
$|K|$ can thus be made to equal any integer greater than
$(p_1|K_1|-1)(p_2|K_2|-1)+|K_3|$. So, if
$\theta_{K}:(p,q,r)\rightarrow S_{|K|}$ is the homomorphism arising from $K$, we have
permutation representations of $(p,q,r)$ for all but finitely many degrees.
By Lemma \ref{cycle.structure} the permutation $\theta_K(x^{-1}y)$ contains 
the $q$-cycle $\mu$ and no other cycles of length
divisible by $q$, so some power of $\theta_{K}(x^{-1}y)$ is just $\mu$. Path-connectedness,
Lemma \ref{primitive} and Theorem \ref{theorem.jordan} give
$\theta_K(G)= A_{|K|}$ or $S_{|K|}$, but the generators of $G$ have odd order, so in fact  
$\theta_K(G)=A_{|K|}$.
\end{proof}

So it remains to give the details. For each of the following cases, the diagrams $K_1,K_2$
and $K_3$ are given and parts $1,2$ and $4$ of the proposition are then easily established. 
Part 3 will prove to be somewhat messier.

\begin{figure}
\begin{center}
\BoxedEPSF{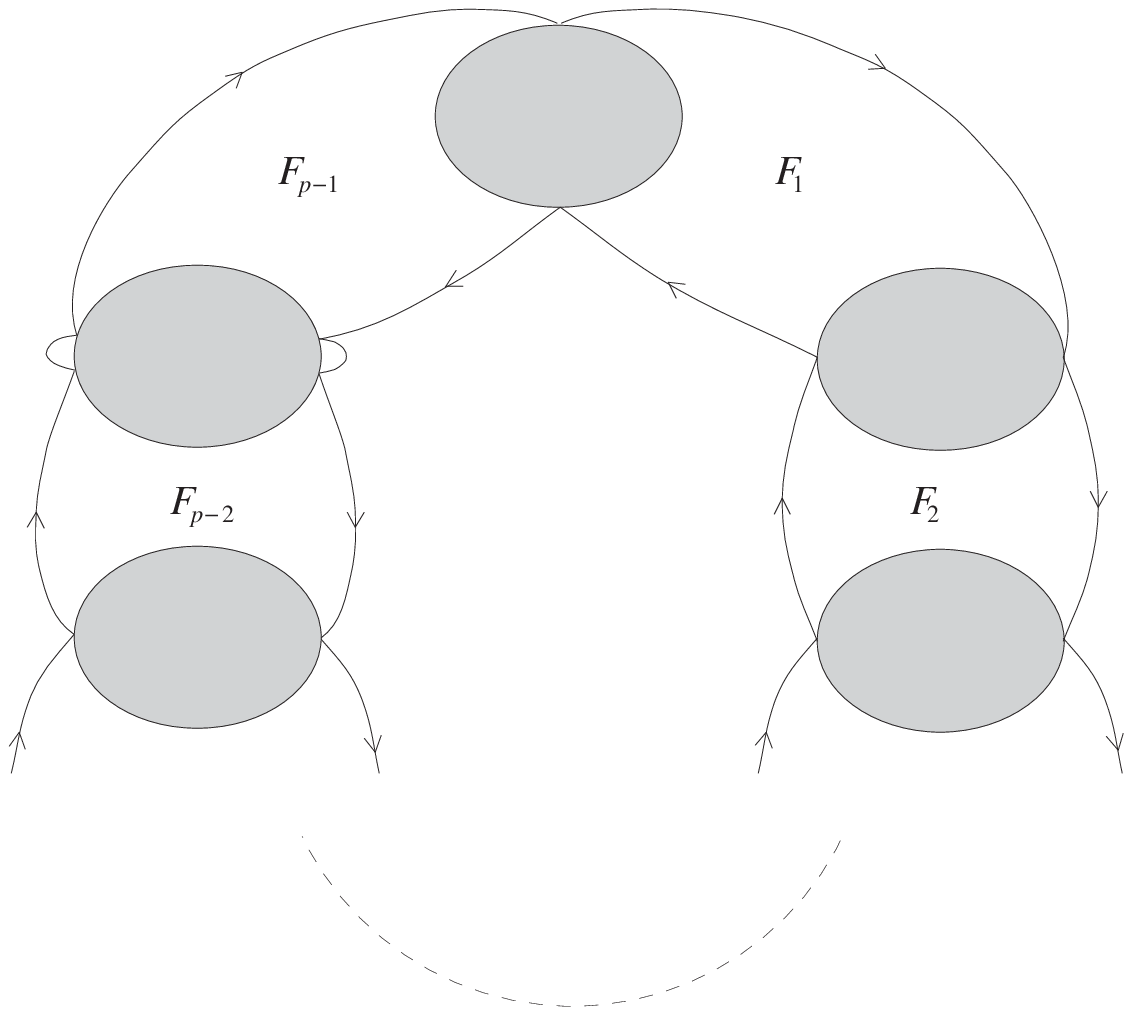 scaled 650}
\caption{} \label{diagramK.1}
\end{center}
\end{figure}

\vspace{1em}

(1). The case $p\geq 7$ and $q\geq p+6$.

\vspace{1em}

Consider Figure \ref{diagramK.1}. We have $q$-gons, $Q_1,\ldots,Q_{p-1}$, with $Q_1$ at
the top and the ordering going clockwise. They are connected by two type $[2,1,\ldots,1]$
$x$-cycles, the number of $1$'s being $p-2$. The connections are such that $Q_i$ contributes one
$y$-arc to the boundary of region $F_{i-1}$, subscripts taken modulo $p-1$. The usual embedding
places Figure \ref{diagramK.1} on the $2$-sphere, as depicted in the picture in fact. The face
$F_{p-1}$ has $q-2$ $y$-arcs in its boundary, faces $F_1,\ldots,F_{p-2}$ have $q$, and there are
four other unshaded faces, two each with label $xy$ and $x^p$.

\begin{figure}
\begin{center}
\BoxedEPSF{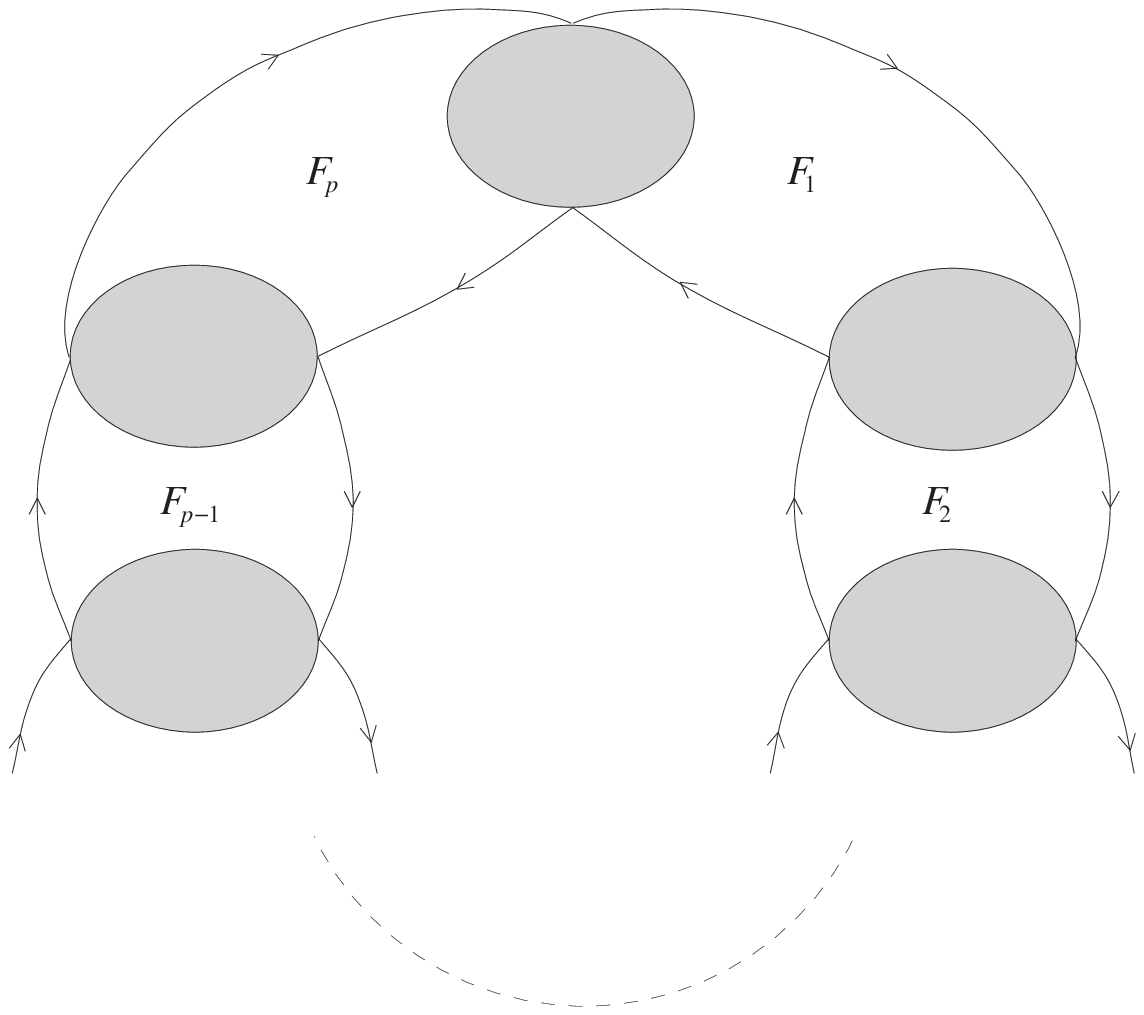 scaled 650}
\caption{} \label{diagramK.3}
\end{center}
\end{figure}

Similarly for Figure \ref{diagramK.3}. We have $q$-gons, $Q_1,\ldots,Q_p$, connected by two type
$[1,\ldots,1]$ $x$-cycles, the number of $1$'s being $p$. The connections are meant to allow $Q_i$
to contribute $\frac{q-1}{2}$ $y$-arcs to the boundary of region $F_{i-1}$, subscripts taken modulo
$p$. The usual embedding places the figure on the 2-sphere also.

Recalling that $q=lp+s$, let $r\geq q+2$ be prime, and $m,\delta$ and $k$ be positive integers
such that,
\begin{itemize}
\item $m$ is largest with $(q+2)+m(p+l+2-s)\leq r$;
\item $\delta$ is largest with $(q+2)+m(p+l+2-s)+\delta (p-3)\leq r$;
\item $k$ is determined by $p-2k+1=r-q-m(p+l+2-s)-\delta (p-3)$.
\end{itemize}
Notice that $2\leq k\leq\frac{p-1}{2}$.
Each $q$-gon $Q_i$ of Figure \ref{diagramK.1} has a number of consecutive free vertices laying
in the boundary of face $F_i$. Assuming for now that this number is sufficient to do so, attach to
$Q_1,\ldots,Q_{p-2}$ arrays of type $\{2^{\delta},k;[2,p-2]^m\}$, and one of type
$\{2^{\delta},k-1;[2,p-2]^m\}$ to $Q_{p-1}$. By (\ref{array.number}) and the definitions of
$m,\delta$ and $k$, each face $F_i$ now has $r$ $y$-arcs in its boundary. We thus have a
spherical $(p,q,r)$-diagram, $K^r_1$. Generally the actual value of $r$ is irrelevant, so we'll
just call this diagram $K_1$.

Take a single $q$-gon, attach to it a type $\{2^{\delta},k;[2,p-2]^m\}$ array and embed. The
resulting spherical $(p,q,r)$-diagram will be our $K_2:=K_2^r$. Our third diagram is slightly more
complicated. In Figure \ref{diagramK.3} attach type $\{k;\mbox{---}\}$ arrays to $Q_1,\ldots,Q_{p-3}$ and
$Q_{p-1}$, using free vertices in the boundary of $F_1,\ldots,F_{p-3}$ and $F_{p-1}$.
To $Q_2,\ldots,Q_{p-2}$ and $Q_p$, attach type $\{2^{\delta};[2,p-2]^m\}$'s, adjacent to
$F_1,\ldots,F_{p-3}$ and $F_{p-1}$, while to $Q_1$ and $Q_{p-2}$, connect
$\{2^{\delta},k;[2,p-2]^m\}$'s adjacent to
$F_p$ and $F_{p-2}$ (the reader should sketch the positions of the various attachments as a guide). 
Again assume for now that there is sufficient space to do all these things. Each
$F_i$ receives $r-q$ new $y$-arcs. The resulting 
$(p,q,r)$-diagram $K_3:=K_3^r$.

Let $N$ be the number of new vertices introduced by an array of type 
$\{2^{\delta},k;[2,p-2]^m\}$. We have
$|K_1|=(p-1)q+(p-2)N+N+1$ and $|K_2|=q+N$. Thus, any common divisor of $|K_1|$ and $|K_2|$ also 
divides 
\begin{equation}\label{rel.prime}
|K_1|-(p-1)|K_2|=1,
\end{equation}
so that $|K_1|$ and $|K_2|$ are relatively prime. Clearly $|K_3|\geq q+3$ and the $K_i$ are 
path-connected.

Let $\s_i$ be as in the proposition, and observe that in $K_3$,
the $\frac{q+1}{2}$ free vertices of $Q_{p-1}$  adjacent to $F_{p-2}$, and the
$\frac{q-1}{2}$ free vertices of $Q_p$ adjacent to $F_p$, form a 
$q$-cycle in $\s_3$. Call any other cycle in $\s_1,\s_2$ or $\s_3$ with length
divisible by $q$ a {\em bad\/} cycle.

We can always arrange things so that bad cycles dissappear and part 3 of the proposition 
thus satisfied. 
The vertices of Figures \ref{diagramK.1}--\ref{diagramK.3} and the $q$-gon
that forms the nucleus of $K_2$ are organised into various cycles. In fact, there are $p-3$
$q$-cycles, a $(q-2)$-cycle and a $(q+2)$-cycle in Figure \ref{diagramK.1}; $p$ $q$-cycles in
Figure \ref{diagramK.3}, and a $q$-cycle in the $q$-gon of $K_2$. A crucial observation is that in
$K_1$ and $K_2$, each of these cycles has exactly one array attached. 
Things are more complicated with $K_3$--one $q$-cycle has
$\{k;\mbox{---}\}$ and $\{2^{\delta},k;[2,p-2]^m\}$ arrays attached, another has 
$\{2^{\delta};[2,p-2]^m\}$ and $\{2^{\delta},k;[2,p-2]^m\}$ arrays, while $p-3$ of them
have $\{k;\mbox{---}\}$ and $\{2^{\delta};[2,p-2]^m\}$. The single $q$-cycle not mentioned
is our precious prime cycle.

We monitor the effect on these cycles of the attached arrays.
First, using the observations
following Lemma \ref{cycle.structure}, one can check that the boosters in a type
$\{k_1,\ldots,k_t;[2,p-2]^m\}$ array contribute bad cycles only when $s=2$. In this
case, the $m$-th booster contains a $q$-cycle. No problem, just modify the chain, and
replacing $X_m$ by $X'_m=[1,p-1]$.

Next the effect of the pendants in an array. Consider one of the $q$-cycles in $K_1$ or $K_2$. If
$m=0$, so that a $\{2^{\delta},k;\mbox{---}\}$ array is attached to the cycle, its length becomes
$$
q-\delta+
\left\{\begin{array}{l}
p-k,\text{ $k$ odd}\\
{\displaystyle -\frac{k}{2}},\text{ $k$ even.}\end{array}\right.
$$
Since $k\leq \frac{p-1}{2}$, we have $p-k\geq \frac{p+1}{2}$, and so the cycle is bad only if
$\delta\geq\frac{p+1}{2}$. The definitions of $m,\delta$ and $k$ give $\delta(p-3)+p-3\leq
p+l+2-s$, so the cycle is bad only if $l-s\geq 7$, that is, $q\geq 8p+1$ (in fact,
$q\geq 2p+1$ will do). By an identical argument, the $(q-2)$-cycle in
$K_1$ becomes bad only if $q\geq 2p+1$, and the $(q+2)$-cycle suffers the same fate under the
addition of a type $\{4;\mbox{---}\}$ array, or only if $q\geq 2p+1$.
Similarly for the $q$-cycles in $K_3$. When $m=0$, we must have $q\geq 3p$ before any turn
bad, and when $m=1$, we must have $q\geq 2p+1$.

What do we do with these bad cycles? When $m\geq 1$ it is simple. Take one of the 
$\{k_1,\ldots,k_t;[2,$ $p-2]^m\}$ arrays attached to the cycle and perform a simultaneous volley
of chain modifications: either replace all $X_i=[2,p-2]$ by $X'_i=[1,p-1]$,
or all $X_i$ by $X'_i=[3,p-3]$, whichever does not create a bad cycle on
the $m$-th booster (they both can't). When $s=2$ and $m\geq 2$, change all $X_i$ to $[1,p-1]$. If
$s=2$ and $m=1$, change $X_1=[1,p-1]$ to $X'_1=[3,p-3]$. In any case the bad cycle
is obliterated and no new bad cycles are created. Remember that when $X'_1=[3,p-3]$, we are assuming
there are two free vertices where the array is attached, but more on this later.

If $m=0$ and a bad cycles arises in $K_3$, spoil one of the attached arrays, assuming for now that
there is enough room to do so. If the bad cycle is in $K_1$ or $K_2$, it would be nice to be rid of
it by spoiling the attached array. Unfortunately, spoiling changes the number of vertices, and
(\ref{rel.prime}) would no longer be valid. So, except for when a $\{4;\text{---}\}$ is attached to
the $(q+2)$-cycle, spoil {\em every\/}
array in these two diagrams (again assuming there is enough room). 
This certainly removes the bad cycle. The danger is that it may have created a new one elsewhere.
If so, remove it by performing a push-pull on the attached array: replace
$\{2^{\delta},k\text{ or }k-1,\frac{p+1}{2};[2,p-2]^m\}$ by
$\{2^{\delta},k-1\text{ or }k-2,\frac{p+3}{2};[2,p-2]^m\}$, or
$\{2^{\delta},1,\frac{p+1}{2};[2,p-2]^m\}$ by $\{2^{\delta},2,\frac{p-1}{2};[2,p-2]^m\}$. In all
the cases that bad cycles arise, $q\geq 2p+1$, so the effect (\ref{push.pull}) of these push-pulls
in both non-trivial and $<q$, so the new bad cycle is removed.

The bad cycle arising when a $\{4;\text{---}\}$ array is attached to the $(q+2)$-cycle in $K_1$ is
removed by similarly spoiling every array in $K_1$ and $K_2$. It can be checked that this creates no
new bad cycles elsewhere. This accounts for all situations where bad cycles arise and establishes
part 3 of the proposition.

Our final task is to see that there are sufficient free vertices in the appropriate places for all
the above to happen. Fix $p$, and for a given $q$, let $\Delta$ be the maximum value
obtained by $\delta$. When $m=0$ the largest number of consecutive free vertices needed anywhere is 
$2\delta+k+\frac{p+1}{2}$: room for a type $\{2^{\delta},k;\text{---}\}$ array and a possible
spoil. Similarly, when $m\geq 1$ we need $2(\delta+1)+k+1$: room for a
$\{2^{\delta},k;[2,p-2]^m\}$ array and a potential  volley of chain
modifications. The $m\geq 1$ needs are less than the $m=0$ needs, and since $k\leq \frac{p-1}{2}$,
these in turn are less than $2\Delta+p$.

Take four consecutive vertices on the $q$-gon of $K_2$ and two on each of
$Q_1,\ldots,$ $Q_{p-2}$ of $K_1$. These are the handles for $K_1$ and $K_2$. Thus, before any arrays
are added, the $q$-gons of $K_1$ and $K_2$ are left with $q-4$ consecutive free vertices.
When $p+6\leq q\leq 2p+1$, we have $\Delta=1$, so $2\Delta+p\leq q-4$, and we are happy.

Now $\Delta$ is the largest multiple of $p-3$ less than $p+l+2-s$.
Thus for a fixed $l$, $\Delta$ and hence $2\Delta+p$ is biggest, and $q-4$ smallest, when $s=1$. It
therefore suffices to show that $2\Delta+p\leq q-4$ for $q=lp+1$. We already have this for $l=2$.
If the inequality is valid for a given $l$, and we increase it by one, then $p+l+2-s$, and hence
$\Delta$, increases by at most one, and so $2\Delta+p$ by at most two. But $q-4$ increases by
$p\geq 7$, and we are home.

In $K_3$ the vertex requirements are greatest and the availability least, on the side of $Q_{p-2}$
adjacent to $F_{p-2}$. 
By considering the possible values of $\Delta$ for $q$ in the range $p+6\leq q\leq 4p-1$, one
can show, using the discussion of when bad cycles arise, that the $\frac{q-3}{2}$ 
consecutive free vertices that are available suffice. For $q\geq 4p+1$, argue as for $K_1$
and $K_2$.

Finally, place a handle on $K_3$ using two vertices of the precious $q$-cycle.

\vspace{1em}

(2). The case $p\geq 7$ and $q=p+2$ or $p+4$.

\vspace{1em}

Diagrams $K_1$ and $K_2$ are the same as in the previous case. That there is sufficient room on $K_1$ and $K_2$ is
a slightly more delicate matter, but the argument is essentially the same. These diagrams can be of no
help to $(11,13,17)$ however, which can be found in \cite[$\S 6$]{Everitt98}.

\begin{figure}
\begin{center}
\BoxedEPSF{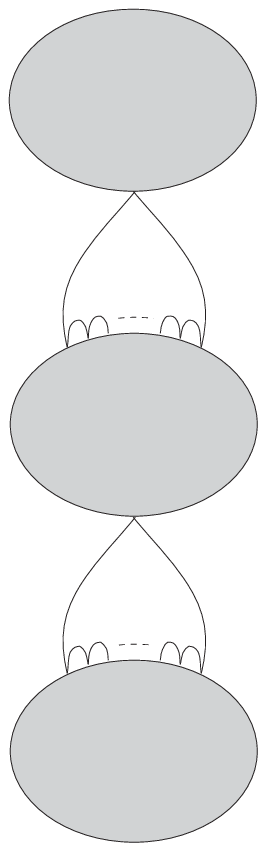 scaled 650}

\caption{} \label{diagramC}
\end{center}
\end{figure}

Unfortunately, there are not enough free vertices on the $K_3$ from case $1$ once $q$ is this close to
$p$. Instead, consider Figure \ref{diagramC}. When $q=p+2$ the large $xy$-face has $q_0=p+10$
$y$-arcs in its boundary, while the minimum $r$ of interest is $r_0=p+6$. For $r>r_0$ prime, let $m$
and $k$ be positive integers such that $m$ is largest with $r_0+m(p+1)\leq r$, and $k$ is determined
by $p-2k+1=r-q_0-m(p+1)$. Add a type $\{k;[1,p-1]^m\}$ array to the  top $q$-gon. The resulting
$(p,q,r)$-diagram will be our $K_3$ for 
$q=p+2$.

Since $k\leq \frac{p+5}{2}$, there is sufficient room on the top $q$-gon for the array with at
least three vertices to spare. Put a handle on the bottom $q$-gon, which also has at least three vertices
to spare.
The middle $q$-gon supplies us with a $q$-cycle. Bad cycles can only arise on the
$q$-gon to which the array is attached. In such a situation, change the two $[1,p-1]$ cycles in 
Figure \ref{diagramC} to type 
$[2,p-2]$'s. This removes the bad cycle.

With $q=p+4$, a $K_3$ diagram for $(7,11,13)$ is in \cite[$\S 6$]{Everitt98}. 
Otherwise the argument is identical with $q_0=p+16$, and
$r_0=p+6$ when $p\geq 13$, or $r_0=17$ when $p=7$. 

\vspace{1em}

(3). The case $p=5$ and $q\geq 17$.

\vspace{1em}

Except for the arrays, diagrams $K_1, K_2$ and $K_3$ are the same as in case 1. For $r\geq q+2$
prime, let $m$ be largest with $(q+2)+m(5+l+2-s)\leq r$; $\delta_1$ largest with 
$(q+2)+m(5+l+2-s)+4\delta_1\leq r$; $\delta_2$ largest with
$(q+2)+m(5+l+2-s)+4\delta_1+2\delta_2\leq r$; and
$k$ determined by $p+2k-1=r-q-m(5+l+2-s)-4\delta_1-2\delta_2$. Add arrays in the same places as
case 1, except replace each $2^{\delta}$ in an array there by $1^{\delta_1},2^{\delta_2}$.
The remainder of the argument is the same.

\vspace{1em}

(4). The case $p=5$ and $q=11,13$.

\vspace{1em}

Diagrams $K_1$ and $K_2$ are as in case 3. For $K_3$, let $r_0=13$ and $q_0=15$ when $q=11$, or
$r_0=17$ and $q_0=21$ when $q=13$. Given $r\geq r_0$ prime, take $m$ largest with $r_0+m(9-s)\leq r$,
and $k$ determined by $p-2k+1=r-q_0-m(9-s)$. Add a type $\{5,k;[1,p-1]^m\}$ array to the top $q$-gon
of Figure \ref{diagramC}, and type $\{5;\text{---}\}$'s to the other two. The resulting $(5,q,r)$-diagram
is our $K_3$. Proceed as in case 2.

\vspace{1em}

(5). The case $p=5$ and $q=7$.

\vspace{1em}

We do $(5,7,11)$ and $(5,7,13)$ in \cite[$\S 6$]{Everitt98}. Diagrams $K_1$ and $K_2$ are the same as in case
1, bar the arrays. Instead, for $r\geq 17$ prime, take $m$ largest with $9+6m\leq r$; $\delta$
largest with $9+6m+2\delta\leq r$, and $k$ given by $5+2k-1=r-7-6m-2\delta$. Somewhat unusually, add
type 
$\{k;[1,4]^{m-1},[1,\bar{\delta},4-\delta]\}$'s and a single type $\{k-1;[1,4]^{m-1},
[1,\bar{\delta},4-\delta]\}$ in all the usual places. When $m\geq 2$ and $\delta=2$, a bad cycle arises 
in the chain of boosters. Remove it by modifying, $X'_m$ being $[2,\bar{2},1]$, and by
replacing the type $2$ pendant on each of the last two boosters by types $1$ and $3$.
For $K_3$, follow the construction of case 2.

\vspace{1em}

(6). The case $p=3$ and $q\geq 17$.

\vspace{1em}

Use figure \ref{diagramK.3}, and allow $Q_i$ to contribute a single $y$-arc to region
$F_{i-1}$. For $r\geq q$ prime, take $m\geq 0$ largest with $q+m(5+l-s)\leq r$ and
$\delta\geq 0$ largest with  $q+m(5+l-s)+2\delta\leq r$. Add type $\{1^{\delta};[2,1]^m\}$
arrays  to each $Q_i$, using the free vertices adjacent to region $F_i$. The resulting
$(3,q,r)$-diagram is our $K_1$.

Spoil the array on $Q_1$, that is, replace by one of type $\{1^{\delta},2;[2,1]^m\}$. This gives
another $(3,q,r)$-diagram, $K_2$. Notice that $|K_1|-|K_2|=1$, so $|K_1|$ and $|K_2|$ are
relatively prime. Place a handle on $Q_2$ and $Q_3$ in each diagram. We can remove bad cycles from
the chains of boosters by the methods of case 1. It is easy to show that none arise elsewhere
in $K_1$. A bad cycle will arise on $Q_1$ in $K_2$ precisely when $m\geq 1$ and $\delta=1$, but
the replacement 
$$\{1^{\delta},2;[2,1]^m\}\rightarrow\{1^{\delta},3;[1,1,\bar{1}],[2,1]^{m-1}\}$$ 
removes it.
The argument of case 1 shows that there are sufficient free vertices for all the arrays and
subsequent modifications.

Take Figure \ref{diagramK.3} with the connecting type $[1,1,1]$ $x$-cycles allowing $Q_i$ to
contribute $\frac{q-1}{2}$ to $F_{i-1}$. Attach type $\{1^{\delta};[2,1]^m\}$ arrays to
$Q_1$ adjacent to $F_1$ and $F_3$, and also to $Q_2$ adjacent to $F_2$. The result is $K_3$.
By the usual argument, there is sufficient room for the arrays as well as to spoil any
array incident with a bad cycle. A $q$-cycle occupies the untouched vertices of $Q_3$ adjacent
to $F_3$ and $Q_2$ adjacent to $F_1$, and a handle for $K_3$ can be safely placed here.   

\vspace{1em}

(7). The case $p=3$ and $q=13$.

\vspace{1em}

You can find $(3,13,17)$ and $(3,13,19)$ in \cite[$\S 6$]{Everitt98}. For $r\geq 23$ prime, use the $K_1$
and $K_2$ of case 6. For $K_3$ attach $\{3^3;\text{---}\}$ arrays to the bottom two $q$-gons
of Figure \ref{diagramC}, and place a handle on the bottom one as well. Place a type
$\{3;\text{---}\}$ on the top $q$-gon. In addition, we need a type $\{1^{\delta};[2,1]^m\}$ on
the top $q$-gon, with $\delta$ and $m$ chosen as in case 6, and this can be spoiled if
necessary to remove bad cycles.

\vspace{1em}

(8). The case $p=3$ and $q=11$.

\vspace{1em}

We do $(3,11,13)$ in \cite[$\S 6$]{Everitt98}. For $r\geq 17$ prime, diagrams $K_1$ and $K_2$ are as in
case 6. For $K_3$ attach type $\{3^2;\text{---}\}$ arrays to the top two $q$-gons
in Figure \ref{diagramC}, and a $\{3^3;\text{---}\}$ array to the bottom. Place a handle on the
middle $q$-gon (which contains our $q$-cycle) and a type $\{1^{\delta};[2,1]^m\}$ array on the top
one. Chose $m$ and $\delta$ according to the usual scheme. Spoil the array to remove any bad
cycles.

\vspace{1em}

(9). The case $p=3$ and $q=7$.

\vspace{1em}

Look in \cite[$\S 6$]{Everitt98} for $(3,7,11)$. For $r\geq 13$ prime, variations on Figure \ref{diagramC}
yield all three diagrams. For consider just the top two $q$-gons and the type $[1,2]$
$x$-cycle connecting them. Place a type $\{1^{\delta};[2,1]^m\}$
array on the top one as usual and handle on each of the top two. The resulting $(3,7,r)$-diagram is
$K_1$. Attach a type
$\{2;\text{---}\}$ array to the bottom $q$-gon. The result is $K_2$. For $K_3$, start from scratch
with  Figure \ref{diagramC}, and attach to the bottom two $q$-gons arrays of type
$\{3;\text{---}\}$, while to the top, attach a type $\{1^{\delta},3;[1,2]^m\}$. Place a handle on 
the bottom $q$-gon.

\vspace{1em}

(10). The case $p=3$ and $q=5$.

\vspace{1em}

You can find $(3,5,7)$ and $(3,5,11)$ in \cite[$\S 6$]{Everitt98}. Otherwise, for
$K_1$ take Figure \ref{diagramC} with type $\{\text{---};[2,1]^m\}$ and
$\{1^{\delta};\text{---}\}$ arrays attached to the second and third $q$-gons respectively,
and with two handles on the top. For $K_2$, place type $\{2;\text{---}\}$ and
$\{1^{\delta};[1,2]^m\}$ arrays on the second and third $q$-gons instead. To get $K_3$, attach a
$\{1^{\delta};[2,1]^m\}$ to  the top $q$-gon and a handle on the bottom one. 

\vspace{1em}

This completes the proof of the theorem.


\begin{acknowledgments}
The author has benefitted from conversations with various people, notably Marston Conder,
Colin Maclachlan, Alan
Reid and Paul Turner. Most of all, I must record a debt of gratitude to Graham Higman, who
provided encouragement and copious improvements to an earlier version of this paper. 
I would also like to thank the referee.
\end{acknowledgments}

\end{article}
\end{document}